\documentclass[11pt]{article}

\usepackage{fullwidth}
\usepackage{xypic}
\usepackage{xy}
\usepackage[top=1.5in, bottom=1.5in, left=1in, right=1in]{geometry}
\usepackage{amsgen}
\usepackage{amsmath}
\usepackage{amstext}
\usepackage{amsbsy}
\usepackage{amsopn}
\usepackage{amsfonts}
\usepackage{amssymb}
\usepackage{eepic}
\usepackage{graphicx}
\usepackage{epsf}
\usepackage{pstricks}

\xyoption{all}

\def\Box{\square}

\def\tra#1{\smash{\mathop{\mid\kern
-1pt\joinrel\relbar\joinrel\relbar}\limits^{*}_{#1}}}
\def\longtra#1{\smash{\mathop{\mid\kern
-1pt\joinrel\relbar\joinrel\relbar\joinrel\relbar}\limits^{*}_{#1}}}
\def\vlongtra#1{\smash{\mathop{\mid\kern
-1pt\joinrel\relbar\joinrel\relbar\joinrel\relbar\joinrel\relbar}\limits^{*}_{#1}}}
\def\vvlongtra#1{\smash{\mathop{\mid\kern
-1pt\joinrel\relbar\joinrel\relbar\joinrel\relbar\joinrel\relbar\joinrel\relbar}\limits^{*}_{#1}}}
\def\vvvlongtra#1{\smash{\mathop{\mid\kern
-1pt\joinrel\relbar\joinrel\relbar\joinrel\relbar\joinrel\relbar\joinrel\relbar\joinrel\relbar}\limits^{*}_{#1}}}
\def\etra#1{\smash{\mathop{\mid\kern
-1pt\joinrel\relbar\joinrel\relbar}\limits_{#1}}}

\def\A{{\cal{A}}}

\def\iff{\Leftrightarrow}
\def\Rw{\Rightarrow}

\def\wt{\widetilde}
\def\wh{\widehat}

\def\M{{\bf M}}
\def\N{\mathbb{N}}

\def\ab{{\bf Ab}}

\def\S{{\cal{S}}}

\def\cl{\mbox{Cl}}

\def\P{{\mathbb{P}}}

\def\G{{\mathbf G}}

\def\V{{\bf V}}

\def\Z{\mathbb{Z}}

\def\p{\varphi}

\def\inv{^{-1}}
\def\la{\langle}
\def\ra{\rangle}

\def\bi{\begin{itemize}}
\def\ei{\end{itemize}}
\def\beq{\begin{equation}}
\def\eeq{\end{equation}}

\def\ds{\displaystyle}
\def\xr{\xrightarrow}

\newtheorem{T}{Theorem}[section]
\newcommand{\bt}{\begin{T}}
\newcommand{\et}{\end{T}}
\newcommand{\ftd}{$\square$\end{T}}

\newtheorem{Proposition}[T]{Proposition}
\newcommand{\bp}{\begin{Proposition}}
\newcommand{\ep}{\end{Proposition}}
\newcommand{\fpd}{$\square$\end{Proposition}}

\newtheorem{Lemma}[T]{Lemma}
\newcommand{\bl}{\begin{Lemma}}
\newcommand{\el}{\end{Lemma}}
\newcommand{\fld}{$\square$\end{Lemma}}

\newtheorem{Corol}[T]{Corollary}
\newcommand{\bc}{\begin{Corol}}
\newcommand{\ec}{\end{Corol}}
\newcommand{\fcd}{$\square$\end{Corol}}

\newtheorem{Result}[T]{Result}
\newcommand{\br}{\begin{Result}}
\newcommand{\er}{\end{Result}}
\newcommand{\frd}{$\square$\end{Result}}

\newtheorem{Example}[T]{Example}
\newcommand{\be}{\begin{Example}}
\newcommand{\ee}{\end{Example}}

\newtheorem{Problem}[T]{Problem}
\newcommand{\bq}{\begin{Problem}}
\newcommand{\eq}{\end{Problem}}

\newtheorem{Remark}[T]{Remark}
\newcommand{\brem}{\begin{Remark}}
\newcommand{\erem}{\end{Remark}}

\newcommand{\proof}
   {\par\medbreak\noindent{\bf Proof}.\enspace}

\newcommand{\qed}{
$\Box$
\par\bigbreak}

\title{The pro-$k$-solvable topology on a free group}\author{{\bf Claude Marion, Pedro V. Silva, Gareth Tracey}}

\date{\today}

\begin{document}

\maketitle

\begin{center}\small

\bigskip

\end{center}

\abstract{We prove that, given a finitely generated subgroup $H$ of a free group $F$, the following questions are decidable: is $H$ closed (dense) in $F$ for the pro-(met)abelian topology? is the closure of $H$ in $F$ for the pro-(met)abelian topology finitely generated? 
We show also that if the latter question has a positive answer, then we can effectively construct a basis for the closure, and the closure has decidable membership problem in any case. Moreover, it is decidable whether $H$ is closed for the pro-$\V$ topology when $\V$ is an equational pseudovariety of finite groups, such as the pseudovariety ${\bf S}_k$ of all finite solvable groups with derived length $\leq k$.
We also connect the pro-abelian topology with the topologies defined by abelian groups of bounded exponent.}

\section{Introduction}
\label{intr}

The classification of finite groups is usually interpreted as the classification of finite simple groups, but universal algebra provides an alternative approach through the concept of a pseudovariety (a class of finite algebras closed under taking subalgebras, homomorphic images and finitary direct products). This is common practice with more general classes of algebras such as semigroups or monoids, where the role played by simple groups has no equivalent.

One of the interesting features of pseudovarieties of finite groups is that they induce a (metrizable) topology on any group $G$. Given such a pseudovariety $\V$, the pro-$\V$ topology on $G$ is the initial topology with respect to all homomorphisms 
$G \to H \in \V$, where $H$ is endowed with the discrete topology. If ${\bf G}$ denotes the pseudovariety of all finite groups, the pro-${\bf G}$ topology is known as the profinite topology.

The profinite topology was introduced by Marshall Hall in \cite{Hal} and he proved in \cite[Theorem 5.1]{Hal2} that every finitely generated subgroup of a free group is closed for the profinite topology.
Over the years, other pseudovarieties $\V$ were considered  \cite{AS,MSW,MSTu,MSTc,MSTs,RZ}, and the following decidability questions became objects of study for an arbitrary finitely generated subgroup $H$ of a free group $F$:
\bi
\item
Can we decide whether $H$ is closed for the pro-$\V$ topology?
\item
Can we decide whether $H$ is dense for the pro-$\V$ topology?
\item
Does the pro-$\V$ closure $\mathrm{Cl}_{\V}(H)$ of $H$ have decidable membership problem?
\item
Can we decide whether $\mathrm{Cl}_{\V}(H)$ is finitely generated and can we compute a basis in the affirmative case?
\ei
To be precise, one should assume that a basis $A$ of $F$ is fixed and we are given as input a finite generating set of $H$, expressed as reduced words over $A \cup A\inv$.

In \cite{RZ}, Ribes and Zalesski\v\i \, answered all these questions positively for the pseudovariety ${\bf G}_p$ of all finite $p$-groups, for an arbitrary prime $p$. In \cite{MSW}, Margolis, Sapir and Weil dealt successfully with the pseudovariety {\bf N} of all finite nilpotent groups (the case of free groups of infinite rank follows from \cite[Corollary 2.4]{MSTc}).
In both cases, the closure of a finitely generated subgroup is always finitely generated.
It is understood that the case of extension-closed pseudovarieties is in general more favourable. 
Indeed, in this case $\mathrm{Cl}_{\V}(H)$ is always finitely generated and of rank at most that of
$H$ \cite[Proposition 3.4]{RZ}. Moreover, all the decidability problems mentioned above are equivalent. This is clearly implicit in \cite{MSW}, but we provide an explicit proof
at the end of Section 2.

Note that ${\bf G}_p$ is extension-closed but {\bf N} is not. And the most famous open problems concern the extension-closed pseudovariety {\bf S} of finite solvable groups.

In this paper, we answer positively all the mentioned problems for
the pseudovariety $\ab$ of finite abelian groups, 
the pseudovariety $\ab(m)$ of finite abelian groups whose exponent divides 
a given positive integer $m$
and the pseudovariety $\M$ of finite metabelian groups, for free groups of arbitrary rank. Obviously, $\ab$ and $\ab(m)$ are not extension-closed, and neither is $\M$: it is easy to check that $S_4 \notin \M$, but $[S_4,S_4] = A_4 \in \M$ and $S_4/A_4 \cong C_2 \in \M$. However, $\ab$, $\ab(m)$ and $\M$ share another favourable property: 
they are equational pseudovarieties of finite groups. In Section 3, we prove some general results concerning such pseudovarieties. In particular, we show that if $\V$ is a decidable equational pseudovariety, then it is decidable whether or not 
$H \leq_{f.g.} F$ is $\V$-closed, where $F$ is any free group.

In Section 4, we apply these results to the pseudovariety ${\bf S}_k$ of all finite solvable groups with derived length $\leq k$ (i.e. finite $k$-solvable groups) for $k \geq 1$. Combining these results with generalizations of theorems by Delgado \cite{Del} and Coulbois \cite{Cou}, we perform in Subsections 4.1 and 4.2 a complete study of the pro-abelian and the pro-metabelian topologies on any free group.

We remark that Delgado's results on the pro-{\bf Ab} topology on free abelian groups were extended by Steinberg in 
\cite{Ste} to any decidable pseudovariety of abelian groups residually containing the integers (namely computing closures of arbitrary rational subsets).

In Section 5, we consider the pro-$\ab(m)$ topology, obtaining similar results to those in Section 4. We also
relate the pro-abelian topology with the pro-$\ab(p^k)$ topologies, where $p$ is a prime and $k$ a positive integer.

\section{Preliminaries}

In this paper, $\P$ denotes the set of all primes. Given a group $G$ and $g,h \in G$, we use the notation 
$[g,h] = g\inv h\inv gh = g\inv g^h$ for the commutator of $g$ and $h$, and denote the derived subgroup of $G$ by $G'$. The second derived subgroup $(G')'$ of $G$ is denoted by $G''$. 
For arbitrary $k \geq 1$, we denote by $G^{(k)}$ the $k$-{\em th derived subgroup} of $G$. A group $G$ is {\em metabelian} if $G'' = 1$. If $G^{(k)} = 1$, then $G$ is $k$-{\em solvable}.

Given a subset $X$ of $G$, we denote by $\la\la X\ra\ra_G$ the normal closure of $X$ in $G$, that is the smallest normal subgroup of $G$ containing $X$.

If $A$ is a set, we set $F_A$ to be the free group on $A$. Also if $H$ is a finitely generated subgroup of $F_A$, we write $H\leq_{f.g.} F_A$.

Let $F_n$ denote the free group of rank $n \in \N$ and fix a basis $A$ of $F_n$. Write $\wt{A} = A \cup A\inv$. We can associate to every finitely generated $H \leq F_n$ a finite automaton $\A(H)$ known as the {\em Stallings automaton} of $H$ \cite{Sta}. For the basic properties of Stallings automata, the reader is referred to \cite[Section 3]{BS}, but we remark that the construction of $\A(H)$ from a finite generating set of $H$ is algorithmically efficient and provides excellent algorithms for the membership problem of $H$ \cite[Proposition 3.5]{BS} and for computing a basis of $H$ \cite[Proposition 3.6]{BS}.

If $H \leq F_n$ is an arbitrary subgroup, we can define the {\em Schreier automaton} $\S(H)$ as follows:
\bi
\item
the vertices are the cosets $Hu$ $(u \in F_n)$;
\item
the edges are of the form $Hu \xr{a} Hua$ $(u \in F_n,\; a \in \wt{A})$;
\item
the vertex 1 is the basepoint (i.e. both initial and terminal vertex).
\ei
Assume that $V$ is the set of all vertices occurring in successful paths with reduced label. If $H$ is finitely generated and we remove from $\S(H)$ all the vertices which are not in $V$, we get the Stallings automaton $\A(H)$. If $H$ is not finitely generated, then $\S(H)$ keeps many of the interesting properties of Stallings automata, but the algorithmic part is usually lost. However, we have the following:

\bp
\label{dmp}
Let $H \leq F_n$ have decidable membership problem. Then we can construct a recursively enumerable basis $B = \{ b_1, b_2, \ldots \}$ of $H$ and express any $h \in H$ as a word on $B$.
\ep

\proof
Let $E_+$ denote the set of positive edges of $\S(H)$ (i.e. with label in $A$). Suppose that $T \subseteq E_+$ is a spanning tree of $\S(H)$ (if we view it as an undirected graph). For every $u \in F_n$, let $\alpha_{Hu}$ label the (unique) $T$-geodesic of the form $1 \xr{\alpha_{Hu}} Hu$. Define
$$B = \{ \alpha_{Hu}a\alpha_{Hua}\inv \mid (Hu,a,Hua) \in E_+ \setminus T\}.$$
The same proof used for Stallings automata \cite[Proposition 3.6]{BS} shows that $B$ is a basis of $H$. We show next how we can use the decidability of the membership problem of $H$ to make the construction of $B$ algorithmic.

Consider the geodesic distance in the Schreier automaton $\S(H)$ and let $D_m(1)$ denote the closed ball with center 1 and radius $m \in \N$. Let $\S_m(H)$ be the subautomaton of $\S(H)$ induced by $D_m(1)$. Clearly, $\S_m(H)$ has less than $\ds\sum_{i=0}^m (2n)^i$ vertices and $\ds\sum_{i=0}^m (2n)^{i+1}$ edges. Since $H$ has decidable membership problem, we can effectively construct $\S_m(H)$ for each $m \in \N$.

Now we can construct a spanning tree $T_m$ of $\S_m(H)$ for each $m \in \N$ with $T_0 \subset T_1 \subset T_2 \subset  \ldots$ It is immediate that $T = \ds\cup_{m\geq 0} T_m$ is a spanning tree for $\S(H)$. Let
$$B_m = \{ \alpha_{Hu}a\alpha_{Hua}\inv \mid Hu,Hua \in D_m(1),\, (Hu,a,Hua) \in E_+ \setminus T_m\}.$$
Then $B = \ds\cup_{m\geq 0} B_m$ is a basis of $H$ and we can enumerate the elements of $B$ as follows: we start by enumerating the elements of $B_0$, then those of $B_1 \setminus B_0$, then those of $B_2 \setminus B_1$ and so on. Since we can compute each $B_m$, we can write $B$ as a sequence $(b_1, b_2,\ldots)$ where each $b_k$ is effectively computable. 

Since every $h \in H$ of length $m$ labels a closed path in $\S_m(H)$ at the basepoint, we can now write it as a word on $B$.
\qed

The following theorem of Karrass and Solitar proves to be useful for our purposes:

\bt
\label{KS}
{\rm \cite[Theorem 1]{KS}}
Let $H$ be a finitely generated subgroup of a free group $F$. Suppose that $H$ contains a nontrivial normal subgroup of $F$. Then $H$ is of finite index in $F$.
\et

In particular, if $F$ is of infinite rank and $H \leq F$ contains a nontrivial normal subgroup of $F$ then $H$ is not finitely generated. Also if $F$ is of finite rank and $H$ contains a nontrivial normal subgroup of $F$, then $H$ is finitely generated if and only if $H$ is of finite index.

A {\em pseudovariety of finite groups} is a class of finite groups closed under taking subgroups, homomorphic images and finitary direct products. Since no other pseudovarieties occur in the paper, we will just write {\em pseudovariety}. 
We consider in this paper the following pseudovarieties:
\bi
\item $\G$ is the pseudovariety of all finite groups;
\item $\ab$ is the pseudovariety of all finite abelian groups;
\item $\M$ is the pseudovariety of all finite metabelian groups;
\item ${\bf S}_k$ is the pseudovariety of all finite $k$-solvable groups (for $k \geq 1$);
\item $\ab(m)$ is the pseudovariety of all finite abelian groups of exponent dividing $m$ (for $m \geq 1$).
\ei
Given a finite group $G$, we denote by $\la G \ra$ the pseudovariety generated by $G$. It contains all the homomorphic images of subgroups of direct powers of the form $G^n$ $(n \geq 1)$. If $m$ is a positive integer, then $\ab(m) = \la C_m \ra$, where $C_m$ denotes the cyclic group of order $m$.

A pseudovariety is {\em decidable} if it has decidable membership problem. For the general theory of pseudovarieties, the reader is referred to \cite{RS}.

Given a pseudovariety $\mathbf{V}$, where we consider finite groups endowed with the discrete topology, the pro-$\mathbf{V}$ topology on a group $G$ is defined as the coarsest topology which makes all morphisms from $G$ into elements of $\mathbf{V}$ continuous. Equivalently, $G$ is a topological group where the normal subgroups $K$ of $G$ such that $G/K\in \mathbf{V}$ form a basis of neighbourhoods of the identity. 

For a topological property $\mathcal{P}$ and a subset $S$ of $F$, we say that $S$ is $\mathbf{V}$-$\mathcal{P}$ if $S$ has property $\mathcal{P}$ in the pro-$\mathbf{V}$ topology on $F$. 

As mentioned in the Introduction, Marshall Hall proved the following seminal result:

\bt
\label{mht}
{\rm \cite[Theorem 5.1]{Hal3}}
Every finitely generated subgroup of a free group is closed for the profinite topology.
\et

The group $G$ is said to be:
\bi
\item
{\em residually finite} if $\{1\}$ is a $\G$-closed subgroup of $G$;
\item
{\em residually} $\V$ if $\{1\}$ is a $\V$-closed subgroup of $G$;
\item
{\em LERF} (locally extended residually finite) if every finitely generated subgroup of $G$ is $\G$-closed.
\ei
In particular, by Theorem \ref{mht} every free group is LERF.

Given a subgroup $H$ of a group $G$, the core $\textrm{Core}_G(H)$ of $H$ in $G$ is the largest normal subgroup of $G$ contained in $H$ and is equal to $\bigcap_{g\in G} g^{-1}Hg$.

The next two fundamental results were proved in full generality by Margolis, Sapir and Weil in \cite{MSW}. Indeed, Marshall Hall's results \cite[Theorems 3.1 and 3.3]{Hal} require $G$ to be residually $\V$. On the other hand, Hall's results cover topologies which do not arise from pseudovarieties, since infinite index subgroups are allowed to be open.

\bt
\label{sopen}
{\rm \cite[Proposition 1.2]{MSW}}
Let $\V$ be a pseudovariety of finite groups. Given a group $G$ and $H \leq G$, the following conditions are equivalent:
\bi
\item[(i)] $H$ is $\mathbf{V}$-open;
\item[(ii)] $H$ is $\mathbf{V}$-clopen;
\item[(iii)] $G/\mathrm{Core}_G(H) \in \mathbf{V}$.
\ei
\et

\bt
\label{sclosed}
{\rm \cite[Proposition 1.3]{MSW}}
Let $\V$ be a pseudovariety of finite groups. Given a group $G$ and $H \leq G$, the following conditions are equivalent:
\bi
\item[(i)] $H$ is $\mathbf{V}$-closed;
\item[(ii)] $H$ is an intersection of $\V$-open subgroups.
\ei
\et

In the particular case of finite index subgroups, we have indeed the following equivalences:

\bp
\label{fi}
Let $\V$ be a pseudovariety of finite groups. Let $G$ be a group and let $H \leq G$ have finite index.
Then the following conditions are equivalent:
\bi
\item[(i)]
$H$ is $\V$-closed;
\item[(ii)]
$H$ is $\V$-clopen;
\item[(iii)]
$G/\mathrm{Core}_G(H) \in \V$.
\ei
\ep

\proof
By Theorem \ref{sopen}, we have (ii) $\iff$ (iii), and (ii) $\Rw$ (i) is trivial. It remains to show that (i) $\Rw$ (ii).

Suppose that $H$ is $\V$-closed. By Theorem \ref{sclosed}, $H$ is an intersection of $\V$-open subgroups. But $[G:H] < \infty$, hence there exist only finitely many subgroups of $G$ containing $H$. It follows that $H$ is an intersection of finitely many $\V$-open subgroups and is therefore $\V$-open.
\qed

We note that a subgroup $H$ of $G$ is $\mathbf{V}$-closed if and only if, for every $g \in G\setminus H$, there exists some $\V$-clopen $K \leq G$ such that $H \leq K$ and $g \notin K$.
On the other hand, a subgroup $H$ of $G$ is $\mathbf{V}$-dense if and only if $HN=G$ for every normal subgroup $N$ of $G$ such that $G/N\in \mathbf{V}$.

The following example shows that Theorem \ref{mht} does not hold for arbitrary subgroups:

\be
\label{exas}
Let $A = \{ a,b\}$ be a basis of $F_2$ and let
$$H = \la a, b^ka^kbab^{-k} \; (k \geq 1)\ra \leq F_2.$$
Then $H$ is not closed for the profinite topology.
\ee

Indeed, we start by remarking that $H < F_2$. Suppose that $b \in H$. Then 
$$b \in H_n = \la a, b^ka^kbab^{-k} \; (k = 1,\ldots,n)\ra \leq F_2$$
for some $n \geq 1$. Since $\A(H_n)$ is of the form
$$
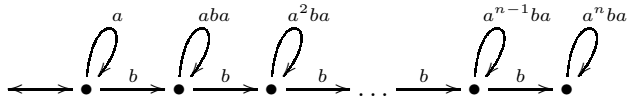
\xymatrix{
\ar@{<->}[r] & \bullet \ar@(u,ur)^a \ar[r]^b & \bullet \ar@(u,ur)^{aba} \ar[r]^b & \bullet \ar@(u,ur)^{a^2ba} \ar[r]^b & \ldots 
\ar[r]^b & \bullet \ar@(u,ur)^{a^{n-1}ba} \ar[r]^b &\bullet \ar@(u,ur)^{a^nba}
}$$
and $\A(H_n)$ does not accept the word $b$, this contradicts \cite[Proposition 3.5]{BS}. Thus $b \notin H$ and so $H < F_2$. Therefore it suffices to show that $H$ is dense for the profinite topology.

Let $N$ be a finite index normal subgroup of $F_2$ with index $n$. Then $a \in H \leq HN$ and $b^n \in N$ yields
$$b = a^{-n}b^{-n}(b^na^nbab^{-n})b^na\inv \in \la H \cup N\ra = HN$$
(since $N \unlhd F_2$). Thus $HN = F_2$ and so $H$ is dense for the profinite topology. Since $H < F_2$, then $H$ is not closed.

\medskip

Given a group $G$ and $S \subseteq G$, we denote by  ${\rm{Cl}}_{\mathbf{V}}^G(S)$ the $\mathbf{V}$-closure of $S$ in $G$. If $G$ is free and no confusion is possible, we will also use the simpler notation  ${\rm{Cl}}_{\mathbf{V}}(S)$. 

We consider now the extension-closed case: 

\bp
\label{t:extclosed}
Let $\V$ be an extension-closed pseudovariety of finite groups. Then the following conditions are equivalent:
\begin{enumerate}
\item[(i)] 
it is decidable whether an arbitrary $H \leq_{f.g.} F_n$ is $\mathbf{V}$-dense;
\item[(ii)]
it is decidable whether an arbitrary $H \leq_{f.g.} F_n$ is $\mathbf{V}$-closed;
\item[(iii)]
${\rm{Cl}}_{\mathbf{V}}(H)$ has a computable (finite) basis for every $H \leq_{f.g.} F_n$.
\end{enumerate}
\ep

\proof
We write $F = F_n$ and fix a basis $A$ of $F$. 
It is clear that condition (iii) implies both  conditions (ii) and (i). Indeed, a subgroup $H$ of $F$ is $\mathbf{V}$-dense if and only if ${\rm{Cl}}_{\mathbf{V}}(H) = F$ and is $\mathbf{V}$-closed if and only if ${\rm{Cl}}_{\mathbf{V}}(H) = H$.

Suppose that condition (ii) holds. We claim that condition (iii) holds. Let $H \leq_{f.g.} F$. Since $\mathbf{V}$ is extension-closed, it follows from \cite[Corollary 2.18]{MSW} that ${\rm{Cl}}_{\mathbf{V}}(H)$ is an {\em overgroup} of $H$ (i.e., its Stallings automaton is a quotient of $\A(H)$ -- obtained through identification of vertices and complete folding). Now $H$ admits only finitely many overgroups (with respect to $A$), which can be effectively computed. By condition (ii), we can compute all the $\mathbf{V}$-closed overgroups of $H$. The smallest such overgroup (equivalently, their intersection) will be ${\rm{Cl}}_{\mathbf{V}}(H)$. Hence condition (iii) holds.

Suppose finally that condition (i) holds. Let $H\leq_{f.g.} F$. 
Let $H_0,H_1,\ldots,H_m$ be the overgroups of $H$, and choose a labelling so that $H_i \not\subseteq H_j$ for all $0 \leq i < j \leq m$.
We show by induction on $i$ that ${\rm{Cl}}_{\mathbf{V}}(H_i) \in \{ H_0,H_1,\ldots,H_i\}$ is computable for $i = 0,\ldots,m$. Since $H_m = H$, the theorem will be proved. 

We consider first the case $i = 0$. Let 
$$B = \{ a \in A \mid a\mbox{ labels some edge in }\mathcal{A}(H)\}.$$ Let $F_B$ be the free group over $B$.
Since $F_B$ is the greatest overgroup of $H$, we have necessarily $H_0 = F_B$. Hence $H_0$ is a free factor of $F$, and note $F$ is trivially $\mathbf{V}$-closed. Since $\mathbf{V}$ is extension-closed, by \cite[Corollary 3.3]{RZ}, $H_0$ is $\mathbf{V}$-closed and so ${\rm{Cl}}_{\mathbf{V}}(H_0) = H_0$.

Assume now that $0 < i \leq m$ and ${\rm{Cl}}_{\mathbf{V}}(H_j)  \in \{ H_0,H_1,\ldots,H_j\}$ is computable for $j = 0,\ldots,i-1$. Let
$$J = \{ j \in \{ 0,\ldots,i-1\} \mid H_i \subset H_j\mbox{ and $H_j$ is $\mathbf{V}$-closed}\}.$$
Note that $0 \in J$.
By Nielsen-Schreier's theorem, each $H_j$ is a free group of finite rank. Thus, by condition (i), we can decide, for each $j \in J$, whether or not $H_i$ is $\mathbf{V}$-dense in $H_j$. Suppose that $H_i$ is $\mathbf{V}$-dense in $H_j$ for some $j \in J$. 
Then $H_j = {\rm{Cl}}_{\mathbf{V}}^{H_j}(H_i)$. But since $\mathbf{V}$ is extension-closed and $H_j$ is $\V$-closed, it follows from \cite[Corollary 3.3]{RZ} that the pro-$\mathbf{V}$ topology on $H_j$ is the subspace topology with respect to the pro-$\mathbf{V}$ topology on $F$. 
Hence $H_j = {\rm{Cl}}_{\mathbf{V}}(H_i)$. 
It follows that $H_i$ is $\mathbf{V}$-dense in $H_j$ for at most one $j \in J$, and in that event we can identify it and deduce that ${\rm{Cl}}_{\mathbf{V}}(H_i) = H_j$. Assume now that there exists no such $j$. Since $\mathbf{V}$ is extension-closed,  ${\rm{Cl}}_{\mathbf{V}}(H_i)$ is an overgroup of $H_i$, and so   $ {\rm{Cl}}_{\mathbf{V}}(H_i) = H_k$ for some $k \in \{ 0,\ldots,m\}$. If $k > i$, then we contradict the rule presiding the enumeration of the overgroups considered. On the other hand, $k < i$ and $H_k$ being $\V$-closed imply that $k \in J$, 
again a contradiction. Thus in this final case we have necessarily $k = i$ and so $ {\rm{Cl}}_{\mathbf{V}}(H_i) = H_i$. Therefore our claim holds for $i$ as required.
\qed

\brem
\label{irec}
Let $\V$ be an extension-closed pseudovariety of finite groups. If the equivalent conditions of Proposition \ref{t:extclosed} hold for every $n \geq 1$, then all conditions hold for arbitrary $H \leq_{f.g.} F$, where $F$ denotes a free group of arbitrary rank.
\erem

\proof
Indeed, let $A$ be a basis of $F$ and let 
$$B = \{ a \in A \mid a \mbox{ or $a\inv$ occurs in the reduced form of some }h \in H\}.$$
Note that $F_B$ contains $H$. Since $H$ is finitely generated, $B$ is finite and $F_B$ is a retract of $F_A$ (i.e. $F_B \leq F_A$ and there is a homomorphism $\theta:F_A \to F_B$ such that $\theta|_{F_B} = 1_{F_B}$). As noted in the proof of \cite[Proposition 2.17]{MSW}, $\V$ being extension-closed implies that $F_B$ is residually $\V$ and so it follows from \cite[Corollary 1.8]{MSW} that $F_B$ is a $\V$-closed subgroup of $F$. 

It follows that if $A$ is infinite, then $H \leq_{f.g.} F$ cannot be $\V$-dense. On the other hand, since $F_B$ is $\V$-closed in $F$ and $\mathbf{V}$ is extension-closed, we know that the pro-$\mathbf{V}$ topology on $F_B$ is the subspace topology with respect to the pro-$\mathbf{V}$ topology on $F$ (see \cite[Corollary 3.3]{RZ}). Thus ${\rm{Cl}}_{\mathbf{V}}^F(H) = {\rm{Cl}}_{\mathbf{V}}^{F_B}(H)$ and 
$H$ is $\V$-closed in $F$ if and only if it is $\V$-closed in $F_B$. This proves our claim.
\qed

\section{Equational pseudovarieties}

In this section, we study the pro-$\V$ topology on a free group $F_A$ for an arbitrary equational pseudovariety $\V$. We recall that $\G$ denotes the pseudovariety of all finite groups.

A {\em variety of groups} is a class of groups closed under taking subgroups, homomorphic images and arbitrary direct products. By Birkhoff's Theorem, a class of groups constitutes a variety if and only if it is the class of all groups satisfying a certain set of group identities. If $X$ denotes a countable alphabet, we say that $u \in F_X$ is a {\em group identity}. A group $G$ satisfies this identity if $\p(u) = 1$ for every homomorphism $\p:F_X \to G$. 
For details on varieties of groups, the reader is referred to \cite{Neu}.

Suppose that $\cal{V}$ is a variety of groups defined by the set of identities $I \subseteq F_X$
We write ${\mathcal{V}} = [I]$.
Given an alphabet $A$, let
$$F_A^{{\cal{V}}} = \la\la \p(u) \mid u \in I,\, \p:F_X \to F_A\mbox{ homomorphism}\ra\ra_{F_A} \unlhd F_A.$$
It follows easily from the definitions and the universal property that $F_A/F_A^{{\cal{V}}}$ is the free object of $\cal{V}$ on $A$. We denote by $\pi_{\cal{V}}:F_A \to F_A/F_A^{{\cal{V}}}$ the canonical homomorphism.

Now ${\cal{V}}^f = {\cal{V}} \cap \G$ is clearly a pseudovariety, known as the {\em finite trace of} $\cal{V}$. In view of Birkhoff's Theorem, a pseudovariety $\V$ is of the form $\mathbf{V} = {\cal{V}}^f$ if and only if it is {\em equational}, that is, there exists a set $I$ of group identities such that $\V$ is the set of all finite groups satisfying all the identities in $I$. Note that arbitrary pseudovarieties require in general pseudoidentities (in view of Reiterman's Theorem). The reader is referred to \cite[Section 7.2]{RS} for further details.

\bt
\label{clvfin}
Let $\cal{V}$ be a variety of groups and let $H \leq_{f.g.} F_A$. Then:
\bi
\item[(i)]
$\mathrm{Cl}_{{\cal{V}}^f}^{F_A}(H) = \pi\inv_{\cal{V}}(\mathrm{Cl}_{\G}^{F_A/F_A^{{\cal{V}}}}(\pi_{\cal{V}}(H))) \geq HF_A^{{\cal{V}}}$;
\item[(ii)]
if $A$ is infinite, then $\mathrm{Cl}_{{\cal{V}}^f}^{F_A}(H)$ is finitely generated if and only if $\cal{V}$ is the variety of all groups;
\item[(iii)]
if $A$ is finite then ${\rm Cl}_{\mathcal{V}^f}^{F_A}(H)$ is finitely generated if and only if $\mathcal{V}$ is the variety of all groups 
or ${\rm Cl}_{\mathcal{V}^f}^{F_A}(H)$ has finite index in $F_A$;
\item[(iv)]
if ${\cal{V}}^f$ is a decidable pseudovariety,
 then it is decidable whether or not $H$ is ${\cal{V}}^f$-closed.
\ei
\et

\proof
Write $\V = {\cal{V}}^f$, $F = F_A$ and $\pi = \pi_{\cal{V}}$.

(i) Let $u \in \mathrm{Cl}_{\V}^{F}(H)$. Assume that $N \unlhd F/F^{{\cal{V}}}$ has finite index and let $K = \pi\inv(N) \unlhd F$. Since $F/K \cong (F/F^{{\cal{V}}})/N \in \V$ and $u \in \mathrm{Cl}_{\V}^{F}(H)$, we get $Ku \cap H \neq \emptyset$. Hence $N\pi(u) \cap \pi(H) \neq \emptyset$ and so $\pi(u) \in \mathrm{Cl}_{\G}^{F/F^{{\cal{V}}}}(\pi(H))$. Thus $u \in \pi\inv(\mathrm{Cl}_{\G}^{F/F^{{\cal{V}}}}(\pi(H)))$ and so $\mathrm{Cl}_{\V}^{F}(H) \leq \pi\inv(\mathrm{Cl}_{\G}^{F/F^{{\cal{V}}}}(\pi(H)))$.

Conversely, let $u \in \pi\inv(\mathrm{Cl}_{\G}^{F/F^{{\cal{V}}}}(\pi(H)))$. Let $N \unlhd F$ be such that $F/N \in \V$. Let $\p:F \to F/N$ be the canonical homomorphism. By the universal property, there exists a homomorphism $\psi:F/F^{{\cal{V}}}\to F/N$ such that $\psi \circ \pi = \p$. If $v \in \pi\inv(\pi(N))$, then 
$$\p(v) = \psi(\pi(v)) \in \psi(\pi(N)) = \p(N) = 1,$$
hence $v \in N$ and so $\pi\inv(\pi(N)) = N$. Now $\pi(N) \unlhd F/F^{{\cal{V}}}$ and 
$$(F/F^{{\cal{V}}})/\pi(N) \cong F/\pi\inv(\pi(N)) = F/N \in \V.$$
Since $\pi(u) \in \mathrm{Cl}_{\G}^{F/F^{{\cal{V}}}}(\pi(H))$, we get $\pi(N) \pi(u) \cap \pi(H) \neq \emptyset$. Hence
there exist some $x \in N$ and $h \in H$ such that $\pi(xu) = \pi(h)$. It follows that $hu\inv \in \pi\inv(\pi(x)) \subseteq \pi\inv(\pi(N)) = N$ and so $Nu \cap H \neq \emptyset$. Thus $u \in \mathrm{Cl}_{\V}^{F}(H)$ and so $\pi\inv(\mathrm{Cl}_{\G}^{F/F^{{\cal{V}}}}(\pi(H))) \leq \mathrm{Cl}_{\V}^{F}(H)$.
This establishes the equality in part (i).

Clearly, $H \leq \mathrm{Cl}_{\V}^{F}(H)$. Since $F^{{\cal{V}}} = \pi\inv(1) \leq \pi\inv(\mathrm{Cl}_{\G}^{F/F^{{\cal{V}}}}(\pi(H)))$, we now get $HF^{{\cal{V}}} \leq \mathrm{Cl}_{\V}^{F}(H)$.

(ii) If $\cal{V}$ is the variety of all groups, then $F^{{\cal{V}}} = \{1\}$ and $\pi$ is the identity, hence
$$\mathrm{Cl}_{\V}^{F}(H) = \pi\inv(\mathrm{Cl}_{\G}^{F/F^{{\cal{V}}}}(\pi(H))) = \mathrm{Cl}_{\G}^{F}(H) = H$$
by part (i) and Theorem \ref{mht}. Therefore $\mathrm{Cl}_{\V}^{F}(H)$ is finitely generated. 

Assume now that $\cal{V}$ is not the variety of all groups. By Birkhoff's Theorem, $\cal{V}$ satisfies a nontrivial identity.  
Thus $F^{{\cal{V}}}$ is nontrivial and it follows from part (i) that $\mathrm{Cl}_{\V}^{F}(H)$ contains a nontrivial normal subgroup of $F$. 
Assume for a contradiction that ${\rm Cl}_{\mathbf{V}}^F(H)$ is finitely generated. By Theorem \ref{KS}, ${\rm Cl}_{\mathbf{V}}^F(H)$ has finite index in $F$. Hence $F$ is finitely generated, a contradiction. Therefore 
$\mathrm{Cl}_{\V}^{F}(H)$ is not finitely generated.

(iii) If $\mathcal{V}$ is the variety of all groups then ${\rm Cl}_{\V}^{F}(H) = H$ by Theorem \ref{mht}, therefore finitely generated. Suppose that $\mathcal{V}$ is not the variety of all groups. Then $F^{\mathcal{V}}$ is nontrivial. Since $F^{{\cal{V}}} \leq \mathrm{Cl}_{\V}^{F}(H)$ by part (i), then
$\mathrm{Cl}_{\V}^{F}(H)$ contains a nontrivial normal subgroup of $F$. Since $F$ has finite rank, then $\mathrm{Cl}_{\V}^{F}(H)$ is finitely generated if and only if $[F:\mathrm{Cl}_{\V}^{F}(H)]$ is finite 
in view of Theorem \ref{KS}. 

(iv) Clearly, $H$ is $\V$-closed if and only if $\mathrm{Cl}_{\V}^F(H) = H$, hence $\mathrm{Cl}_{\V}^F(H)$ being finitely generated is a necessary condition. 

Suppose first that $A$ is infinite. By part (ii), $\mathrm{Cl}_{\V}^F(H)$ being finitely generated can only happen if $\V = \G$, and in that case $H$ is $\G$-closed by Theorem \ref{mht}.

Assume now that $A$ is finite.
In view of part (ii), this can only happen if $[F:H]$ is finite. This necessary condition can easily be checked with the help of the Stallings automaton $\A(H)$ (which must be complete by \cite[Proposition 3.8]{BS}). Thus we may assume that $[F:H]$ is finite.
Now, with the help of $\A(H)$, it is easy to compute the Stallings automaton of $C = \textrm{Core}_F(H)$: we consider the finitely many conjugates of $H$ and intersect  them using the direct product of their Stallings automata. Since $C \unlhd F$, then the underlying graph of $\A(C)$ is indeed the Cayley graph of $F/C$ with respect to the alphabet of $F$, hence $F/C$ is computable. Since $\V$ is a decidable pseudovariety, then we can decide whether or not $F/C \in \V$. By Proposition \ref{fi}, this is equivalent to $H$ being $\V$-closed and we are done.
\qed

The last result of this section will prove useful when considering free objects of arbitrary rank:

\bp
\label{lerfir}
Let $\cal{V}$ be a variety of groups with residually finite free objects. If $F_A/F_A^{{\cal{V}}}$ is LERF for every finite alphabet $A$, then it is LERF for every alphabet.
\ep

\proof
Let $A$ be an arbitrary alphabet and let $\alpha:F_A \to F_A/F_A^{{\cal{V}}}$ be the canonical homomorphism. Let $H \leq_{f.g.} F_A$. As in the discussion of Remark \ref{irec}, let $B$ be the (finite) set of all letters of $A$ occurring in the reduced form of some $u \in H$. We noted then that $F_B$ is a retract of $F_A$, now we claim that $F_B/F_B^{{\cal{V}}}$ can be viewed as a retract of $F_A/F_A^{{\cal{V}}}$. 

Indeed, let $\theta:F_A \to F_B$ be the homomorphism fixing each $b \in B$ and sending each $a \in A\setminus B$ to 1. Let $\beta:F_B \to F_B/F_B^{{\cal{V}}}$ be the canonical homomorphism. Considering the composition $\beta \circ \theta: F_A \to F_B/F_B^{{\cal{V}}} \in {\cal{V}}$, it follows from the universal property that there exists some homomorphism $\p:F_A/F_A^{{\cal{V}}} \to F_B/F_B^{{\cal{V}}}$ such that the diagram
$$\xymatrix{
F_A \ar[rr]^{\beta\circ\theta} \ar[dd]_{\alpha} && F_B/F_B^{{\cal{V}}}\\ &&\\
F_A/F_A^{{\cal{V}}} \ar[uurr]_{\p} &&
}$$
commutes. 
Similarly, by considering the inclusion $\iota:F_B \to F_A$, we show that there exists some homomorphism $\psi:F_B/F_B^{{\cal{V}}} \to F_A/F_A^{{\cal{V}}}$ such that the diagram
$$\xymatrix{
F_B \ar[rr]^{\alpha\circ\iota} \ar[dd]_{\beta} && F_A/F_A^{{\cal{V}}}\\ &&\\
F_B/F_B^{{\cal{V}}} \ar[uurr]_{\psi} &&
}$$
commutes. 

Checking through letters, we get $\p \circ \psi = 1_{F_B/F_B^{{\cal{V}}}}$, hence $\psi$ is injective and so $F_B/F_B^{{\cal{V}}} \cong \psi(F_B/F_B^{{\cal{V}}})$. Moreover, they are homeomorphic when endowed with the profinite topology. 

Now $\beta(H)$ is a finitely generated subgroup of $F_B/F_B^{{\cal{V}}}$. Since $F_B/F_B^{{\cal{V}}}$ is LERF by hypothesis, then $\beta(H)$ is a $\G$-closed subgroup of $F_B/F_B^{{\cal{V}}}$. Hence $\alpha(H) = \alpha(\iota(H)) = \psi(\beta(H))$ is a $\G$-closed subgroup of $\psi(F_B/F_B^{{\cal{V}}})$.

On the other hand, it follows from $\p \circ \psi = 1_{F_B/F_B^{{\cal{V}}}}$ that
$\psi \circ\p \circ \psi = \psi$, thus $\psi(F_B/F_B^{{\cal{V}}})$ is a retract of $F_A/F_A^{{\cal{V}}}$.  
By \cite[Proposition 1.6]{MSW}, the profinite topology on $\psi(F_B/F_B^{{\cal{V}}})$ is the subspace topology with respect to the profinite topology on $F_A/F_A^{{\cal{V}}}$. Thus $\alpha(H) = C \cap \psi(F_B/F_B^{{\cal{V}}})$ for some $\G$-closed subset $C$ of $F_A/F_A^{{\cal{V}}}$.
Since $F_A/F_A^{{\cal{V}}}$ is residually finite by hypothesis, it follows from \cite[Corollary 1.8]{MSW} that $\psi(F_B/F_B^{{\cal{V}}})$ is a $\G$-closed subgroup of $F_A/F_A^{{\cal{V}}}$, hence $\alpha(H)$ is itself a $\G$-closed subgroup of $F_A/F_A^{{\cal{V}}}$. 

Now every finitely generated subgroup of $F_A/F_A^{{\cal{V}}}$ is of the form $\alpha(H)$ for some finitely generated subgroup $H$ of $F_A$ (we can lift through $\alpha$ a finite set of generators).
Therefore $F_A/F_A^{{\cal{V}}}$ is LERF.
\qed

\section{The pro-${\bf S}_k$ topology}

In this section, we study the pro-${\bf S}_k$ topology. Let $\wh{{\bf S}_k}$ denote the variety of all $k$-solvable groups. The following lemma is folklore:

\bl
\label{skep}
For every alphabet $A$ and every $k \geq 1$, we have:
\bi
\item[(i)] ${\bf S}_k = \wh{{\bf S}_k}^f$ and is therefore an equational pseudovariety;
\item[(ii)] $F_A^{\wh{{\bf S}_k}} = F_A^{(k)}$.
\ei
\el

In fact, each $\wh{{\bf S}_k}$ can be defined by a single identity: $[x_1,y_1]$ for $k = 1$, $[[x_1,y_1],[x_2,y_2]]$ for $k = 2$, and so on, and this yields (ii). Moreover, $\mathbf{S}_k$ is a decidable pseudovariety.

Now we can apply Theorem \ref{clvfin} to get:

\bt
\label{prosk}
Let $H \leq_{f.g.} F_A$ and $k \geq 1$. Then:
\bi
\item[(i)]
$\mathrm{Cl}_{{\bf S}_k}(H) \geq HF_A^{(k)}$;
\item[(ii)]
if $A$ is infinite, then $\mathrm{Cl}_{{\bf S}_k}(H)$ is not finitely generated;
\item[(iii)]
if $A$ is finite, then $\mathrm{Cl}_{{\bf S}_k}(H)$ is finitely generated if and only if it has finite index in $F_A$;
\item[(iv)]
it is decidable whether or not $H$ is ${\bf S}_k$-closed.
\ei
\et

We can progress further in the abelian and metabelian cases. A crucial difference is that free abelian and free metabelian groups of finite rank are LERF, and this fails for the free $k$-solvable group $F_n/F_n^{(k)}$ for $ n\geq 2$ and $k \geq 3$ \cite{Aga,Gru}. Moreover, Umirbaev proved that the membership problem for finitely generated subgroups is undecidable for $F_2/F_2^{(3)}$, unlike the abelian and metabelian cases \cite{Umi}.

\subsection{The pro-abelian topology}

We establish next all the basic facts concerning the pro-abelian topology. But first we extend Delgado's theorem to free abelian groups of arbitrary rank:

\bl
\label{fabir}
Free abelian groups are LERF.
\el

\proof
The free abelian group of basis $A$ can be described as the direct sum $\ds\oplus_{a \in A} \Z$ (i.e., all mappings $A \to \Z$ with finite support and componentwise sum). Since $\ds\oplus_{a \in A} \Z$ is clearly residually finite, it follows from \cite{Del} and Proposition \ref{lerfir} that $\ds\oplus_{a \in A} \Z$ is LERF for any alphabet $A$.
\qed

\bt
\label{proab}
Let $H \leq_{f.g.} F_A$. Then:
\bi
\item[(i)]
$\mathrm{Cl}_{\ab}(H) = HF'_A$;
\item[(ii)]
if $A$ is infinite, then $\mathrm{Cl}_{\ab}(H)$ is not finitely generated;
\item[(iii)] 
if $A$ is finite, then $\mathrm{Cl}_{\ab}(H)$ is finitely generated if and only if it has finite index in $F_A$;
\item[(iv)]
it is decidable whether or not $\mathrm{Cl}_{\ab}(H)$ is finitely generated;
\item[(v)]
if $\mathrm{Cl}_{\ab}(H)$ is finitely generated, we can compute a basis of $\mathrm{Cl}_{\ab}(H)$;
\item[(vi)]
it is decidable whether or not $H$ is $\ab$-closed;
\item[(vii)]
it is decidable whether or not $H$ is $\ab$-dense;
\item[(viii)]
the membership problem for $\mathrm{Cl}_{\ab}(H)$ is decidable;
\item[(ix)]
if $A$ is finite, then we can define a recursively enumerable basis $B$ of $\mathrm{Cl}_{\ab}(H)$ and produce an algorithm which expresses any element of $\mathrm{Cl}_{\ab}(H)$ as a word on $B$.
\ei
\et

\proof
Write $F = F_A$ and let $\pi:F \to F/F'$ be the canonical homomorphism. 

(i) By Lemma \ref{skep}(ii), the free abelian group over $A$ is $F/F'$. By Lemma \ref{fabir}, $\pi(H)$ is closed in $F/F'$ for the profinite topology. If we take $\cal{V}$ to be the variety of all abelian groups, it follows from Theorem \ref{clvfin}(i) that
\beq
\label{proab1}
\cl_{\ab}(H) = \mathrm{Cl}_{{\cal{V}}^f}^{F}(H) = \pi\inv(\mathrm{Cl}_{\G}^{F/F'}(\pi(H))) = \pi\inv(\pi(H)) = HF'
\eeq
as claimed.

(ii) By Theorem \ref{prosk}(ii).

(iii) By Theorem \ref{prosk}(iii).

(iv) In view of part (ii), we may assume that $A$ is finite.
It follows from (\ref{proab1}) and the homomorphism theorems that
\beq
\label{proab2}
F/\mathrm{Cl}_{\ab}(H) = F/\pi\inv(\pi(H)) \cong (F/F')/\pi(H).
\eeq
We may assume that $F/F' = \mathbb{Z}^n$ for $n = |A|$, hence $\pi(H) \leq \mathbb{Z}^n$.
By \cite[Chapter II, Theorem 1.6]{Hun}, there exist integers $e_1,\dots,e_n \geq 0$ such that $e_1\mid e_2\mid\dots\mid e_n$, and an automorphism $\phi$ of $\mathbb{Z}^n$ such that 
$$\phi(\pi(H))=\bigoplus_{i=1}^n e_i\mathbb{Z}.$$ Moreover one can explicitly produce $e_1,\ldots,e_n$ and $\phi$. 
It is immediate that $\mathbb{Z}^n/\pi(H)$ is finite if and only if $e_n > 0$, hence we can decide whether or not $\mathbb{Z}^n/\pi(H)$ is finite. It follows from (\ref{proab2}) that we can decide whether or not $\mathrm{Cl}_{\ab}(H)$ has finite index in $F_n$. By part (iii), we can decide whether or not $\mathrm{Cl}_{\ab}(H)$ is finitely generated.

(v) If $\mathrm{Cl}_{\ab}(H)$ is finitely generated, then $A$ is finite by part (ii). Write $n = |A|$. Then 
$\mathrm{Cl}_{\ab}(H)$
has finite index in $F$ by part (iii) and it follows from its proof that
we can effectively compute the finite group $\mathbb{Z}^n/\pi(H)$, which is isomorphic to $F/\mathrm{Cl}_{\ab}(H)$ by (\ref{proab2}). 
Hence $\mathrm{Cl}_{\ab}(H)$ is a normal subgroup of finite index of $F$. It follows that the Stallings automaton $\A$ of $\mathrm{Cl}_{\ab}(H)$ with respect to $A$ is precisely the (computable) Cayley graph of $F/\mathrm{Cl}_{\ab}(H) \cong \mathbb{Z}^n/\pi(H)$ with respect to $A$, which can now be used to 
produce a basis of $\mathrm{Cl}_{\ab}(H)$.

(vi) By Theorem \ref{prosk}(iv).

(vii) Suppose first that $A$ is infinite. Suppose that $H$ is $\ab$-dense. Then $\mathrm{Cl}_{\ab}(H) = HF' = F$.  Also $HF'/F'$ is finitely generated, since $H$ is. It follows that $F/F'$ is finitely generated, a contradiction since $F$ is of infinite rank.
Hence we may assume that $A$ is finite. We have that $H$ is $\ab$-dense if and only if $\mathrm{Cl}_{\ab}(H) = F$, and we may assume that $\mathrm{Cl}_{\ab}(H)$ is finitely generated and the finite quotient $F/\mathrm{Cl}_{\ab}(H)$ is computable. Therefore we only need to check whether or not this finite group is trivial.

(viii) Since $\cl_{\ab}(H) = \pi\inv(\pi(H))$ by (\ref{proab1}), it suffices to note that the membership problem is decidable for every (finitely generated) subgroup of $F/F' \cong \ds\oplus_{a \in A} \Z$. This follows from a well-known fact: it is decidable whether or not a system of finitely many linear diophantine equations admits a solution over the integers.

(ix) By part (viii) and Proposition \ref{dmp}.
\qed

\subsection{The pro-metabelian topology}

We present next similar results for the pro-metabelian topology.

We start with a useful lemma. Given subgroups $H,K$ of a group $G$, we denote by $[H,K]$ the subgroup of $G$ generated by all commutators of the form $[h,k]$ with $h \in H$ and $k \in K$.

\bl
\label{derp}
Let $G$ be a group such that $G = NH$ for some $N \unlhd G$ and $H \leq G$. Then $G' = N'[N,H]H'$.
\el
 
 \proof
Let $a,b \in G$ and write $a=nh$ and $b=mk$ for $n, m \in N$ and $h, k \in H$. Then
$$[a,b] = [nh,mk] = [n,mk]^h[h,mk] = ([n,k][n,m]^k)^h([h,k][h,m]^k),$$
using the commutator identities.
 
Now, $[n,kh] = [n,h][n,k]^h$, so
\beq
\label{derp1}
[n,k]^h = [n,h]^{-1}[n,kh] \in [N,H].
\eeq
Also, $[n,m]^{kh} \in N'$ since $N \unlhd G$. Next, $[h,k] \in H'$. 
And finally, $[h,m]^k = ([m,h]^{k})\inv \in [N,H]$ by (\ref{derp1}).
 
Thus, $G' \leq [N,H]N'H'[N,H]$. The argument at (\ref{derp1}) shows that $[N,H]$ is normalized by $H$. It is also normalized by $N$, since $[nm,h]=[n,h]^m[m,h]$ yields $[n,h]^m=[nm,h][m,h]^{-1}$.
 
Therefore $[N,H] \unlhd NH=G$, so $G' \leq [N,H]N'H'[N,H]=N'[N,H]H'$.
 
Since the opposite containment is trivial, the claim follows.
\qed

Next we extend Coulbois's theorem to free metabelian groups of arbitrary rank:

\bl
\label{fmir}
Free metabelian groups are LERF.
\el

\proof
Let $A$ be an arbitrary alphabet and let $G = F_A/F_A^{\wh{\M}}$ be the free metabelian group on $A$ (i.e. the free object of $\wh{\M}$ on $A$). Let $g \in G \setminus \{ 1\}$. By \cite[36.35]{Neu}, $G$ is {\em residually 2-generated}, hence there exists some $N \lhd G$ such that $G/N$ is 2-generated and $g \notin N$. Thus $G/N$ is a finitely generated metabelian group. By a theorem of Philip Hall \cite{Hal2}, $G/N$ must be residually finite. Since $gN \neq N$, there exists some $K \lhd G/N$ such that $(G/N)/K$ is finite and $gN \notin K$.

Let $\pi:G \to G/N$ denote the canonical homomorphism. Then $\pi\inv(K) \lhd G$ and $G/\pi\inv(K) \cong (G/N)/K$ is finite. Since $gN \notin K$ yields $g \notin \pi\inv(K)$, this shows that $G$ is residually finite.

Therefore free metabelian groups are residually finite, and it follows from \cite{Cou} (see also \cite{Ali}) and Proposition \ref{lerfir} that all free metabelian groups are LERF.
\qed

We can now prove our main result:

\bt
\label{prometa}
Let $H \leq_{f.g.} F_A$. Then:
\bi
\item[(i)]
$\mathrm{Cl}_{\M}(H) = HF''_A$;
\item[(ii)] 
if $A$ is infinite, then $\mathrm{Cl}_{\M}(H)$ is not finitely generated;
\item[(iii)] 
if $A$ is finite, then $\mathrm{Cl}_{\M}(H)$ is finitely generated if and only if it has finite index in $F_A$;
\item[(iv)]
it is decidable whether or not $\mathrm{Cl}_{\M}(H)$ is finitely generated;
\item[(v)]
if $\mathrm{Cl}_{\M}(H)$ is finitely generated, we can compute a basis of $\mathrm{Cl}_{\M}(H)$;
\item[(vi)]
it is decidable whether or not $H$ is $\M$-closed;
\item[(vii)]
it is decidable whether or not $H$ is $\M$-dense;
\item[(viii)]
the membership problem for $\mathrm{Cl}_{\M}(H)$ is decidable.
\item[(ix)]
if $A$ is finite, then we can define a recursively enumerable basis $B$ of $\mathrm{Cl}_{\M}(H)$ and produce an algorithm which expresses any element of $\mathrm{Cl}_{\M}(H)$ as a word on $B$.
\ei
\et

\proof
Write $F = F_A$. Let $\pi:F \to F/F''$ be the canonical homomorphism. 

(i) By Lemma \ref{fmir}, $\pi(H)$ is closed in $F/F''$ for the profinite topology. If we take $\cal{V}$ to be the variety of all metabelian groups, it follows from Theorem \ref{clvfin}(i) that
\beq
\label{prometa1}
\cl_{\M}(H) = \mathrm{Cl}_{{\cal{V}}_{{\rm fin}}}^{F}(H) = \pi\inv(\mathrm{Cl}_{\G}^{F/F''}(\pi(H))) = \pi\inv(\pi(H)) = HF''
\eeq
as claimed.

(ii) By Theorem \ref{prosk}(ii).

(iii) By Theorem \ref{prosk}(iii).

(iv) In view of part (ii), we may assume that $A$ is finite.
Let $G = HF' = F'H \leq F$. Since $F' \unlhd G$ and $H \leq G$, it follows from Lemma \ref{derp} that
\beq
\label{metcl1}
G' = F''[F',H]H'.
\eeq

It is easy to check that
\beq
\label{metcl2}
G/(HF'') \cong F'/(F' \cap HF'').
\eeq
Indeed, the kernel of the canonical homomorphism $\p:F' \to G/(HF'')$ is $F' \cap HF''$ and surjectivity follows from $uhHF'' = uHF''$ holding for all $u \in F'$ and $h \in H$. 

Since $F'' \leq F' \cap HF''$, then $F'/(F' \cap HF'')$ is abelian and so is $G/(HF'')$. Hence $G' \leq HF''$ and so $HG' \leq HF''$. Since $F'' \leq G'$ by (\ref{metcl1}), it follows that $HG' = HF''$.
 
By part (i), we have $\mathrm{Cl}_{\M}(H) = HF'' = HG'$. Also $\mathrm{Cl}_{\ab}(H) = HF' = G$ by Theorem \ref{proab}(i). By part (ii), we want to decide whether or not $[F:HG'] < \infty$.
Since $HG' \leq G \leq F$, we have
$$[F:HG'] < \infty \mbox{ if and only if $([F:G] < \infty$ and }[G:HG'] < \infty).$$
By Theorem \ref{proab}, we can decide whether or not $[F:G] < \infty$, and in the affirmative case we can compute a basis for $G$ (which is a free group of finite rank in view of Nielsen's theorem). Thus we may assume that this is the case, and we can write any given element of $H$ in terms of this basis. But now $HG' = \mathrm{Cl}_{\ab}^{G}(H)$ and once again we can decide whether or not $[G:HG'] < \infty$. Therefore we can 
decide whether or not $[F:HG'] < \infty$ as required.

(vi) By Theorem \ref{prosk}(iv).

(v) If $C = \mathrm{Cl}_{\M}(H) = HF''$ is finitely generated, then $A$ is finite 
by part (ii). Write $n = |A|$. Then 
$C$
has finite index in $F$ by part (iii). 
Using the notation $G = HF'$ from the proof of (iv), we have $C = HG'$ and
\beq
\label{prin}
[F:C] = [F:G][G:HG'].
\eeq
Moreover, $F/G$ and $G/HG'$ are both finite abelian groups. Note that $G = \mathrm{Cl}_{\ab}(H)$ and $G$ is finitely generated as a finite index subgroup of $F$. Hence $[F:G]$ is computable in view of Theorem \ref{proab}(v), with the help of Stallings automata: by \cite[Proposition 3.8]{BS}, $[F:G]$ is the number of vertices of $\A(G)$. 

Since $G$ is finitely generated, using Theorem \ref{proab}(ix), we can compute a basis of $G$ and can write any given element of $H$ in terms of this basis.
We can also compute $[G:HG'] = [G:\mathrm{Cl}_{\ab}^G(H)]$. By (\ref{prin}), $[F:C]$ is computable. Now $F$ has only finitely many subgroups $K_1,\ldots,K_s$ of index $m = [F:C]$, defined through their Stallings automata, which must be complete and possess $m$ vertices. Thus it suffices to decode whether each of these $K_i$ equals $C$. We have two obvious necessary conditions: $H \leq K_i$ and $K_i$ being $\M$-closed. The first is decidable through Stallings automata, and the second in view of part (vi). Since $[F:C] = m = [F:K_i]$ for each $i$, these two necessary conditions are also sufficient, hence we can effectively compute a basis of $C$.

(vii) If $A$ is infinite, then we know by the proof of Theorem \ref{proab}(vii) that $H$ is not $\ab$-dense. Since 
$\mathrm{Cl}_{\M}(H) = HF'' \leq HF' = \mathrm{Cl}_{\ab}(H)$, $H$ cannot be $\M$-dense either. 
Hence we may assume that $A$ is finite. We have that $H$ is $\M$-dense if and only if $\mathrm{Cl}_{\M}(H) = F$, hence decidability follows from parts (iv) and (v).

(viii) Suppose first that $A$ is finite.
Since $\cl_{\M}(H) = \pi\inv(\pi(H))$ by (\ref{prometa1}), it suffices to note that the membership problem is decidable for every finitely generated subgroup of the free metabelian group $F/F''$,  a theorem proved by Romanovski\v\i \, \cite{Rom}.

Assume now that $A$ is infinite. Let $u \in F_A/F_A''$.
In view of Lemma \ref{fmir}, free metabelian groups are residually finite and the hypotheses of Proposition \ref{lerfir} are satisfied. Adapting its proof, let $B$ be the (finite) set of all letters of $A$ occurring in the reduced form of some $h \in H$ or $u$. The proof shows that $F_B/F_B''$ can be viewed as a retract of $F_A/F_A''$, hence membership of $u$ can be tested within $F_B/F_B''$, and we are done by the previous case.

(ix) By part (viii) and Proposition \ref{dmp}.
\qed

\section{The pro-$\ab(m)$ topology}

In this section, we adapt the proof of Theorem \ref{proab} to the pro-$\ab(m)$ topology and establish connections with the pro-abelian topology.

\bt
\label{proabm}
Let $m \geq 1$ and $H \leq_{f.g.} F_A$. Then:
\bi
\item[(i)]
$\mathrm{Cl}_{\ab(m)}(H) = HF'_AF_A^m$;
\item[(ii)]
$\mathrm{Cl}_{\ab(m)}(H)$ is finitely generated if and only if $A$ is finite;
\item[(iii)]
it is decidable whether or not $\mathrm{Cl}_{\ab(m)}(H)$ is finitely generated;
\item[(iv)]
if $\mathrm{Cl}_{\ab(m)}(H)$ is finitely generated, we can compute a basis of $\mathrm{Cl}_{\ab(m)}(H)$;
\item[(v)]
it is decidable whether or not $H$ is $\ab(m)$-closed;
\item[(vi)]
$H$ is $\ab(m)$-dense if and only if $H$ is $\ab(p)$-dense for every prime $p$ dividing $m$;
\item[(vii)]
it is decidable whether or not $H$ is $\ab(m)$-dense;
\item[(viii)]
the membership problem for $\mathrm{Cl}_{\ab(m)}(H)$ is decidable;
\item[(ix)]
if $A$ is finite, then we can define a recursively enumerable basis $B$ of $\mathrm{Cl}_{\ab(m)}(H)$ and produce an algorithm which expresses any element of $\mathrm{Cl}_{\ab(m)}(H)$ as a word on $B$.
\ei
\et

\proof
Write $F = F_A$. 

(i) Let ${\cal{V}} = [ x\inv y\inv xy,x^m]$ denote the variety of all abelian groups of exponent dividing $m$. Then 
$$F^{{\cal{V}}} = \la\la \p([x,y]), \p(x^m) \mid \p:F_X \to F\mbox{ homomorphism}\ra\ra_F
= \la\la [u,v],u^m \mid u,v \in F\ra\ra_F.$$
Since $F'$ and $F^m = \la u^m \mid u \in F\ra$ are both normal subgroups of $F$, we get $F^{{\cal{V}}} = F'F^m$.
Hence the free object of $\mathcal{V}$ over $A$ is $F/(F'F^m)$, and it is easy to see that $F/(F'F^m) \cong C_m^A$. Note that $\ab(m) = {\cal{V}}^f$.

If $A$ is finite, then $C_m^A$ is finite and therefore LERF. Since $C_m^A$ is residually finite for arbitrary $A$, it follows from Proposition \ref{lerfir} that $C_m^A$ is always LERF. 
Let $\pi:F \to C_m^A$ be the canonical homomorphism. It follows from Theorem \ref{clvfin}(i) that
\beq
\label{proabm1}
\cl_{\ab(m)}(H) = \mathrm{Cl}_{{\cal{V}}^f}^{F}(H) = \pi\inv(\mathrm{Cl}_{\G}^{C_m^A}(\pi(H))) = \pi\inv(\pi(H)) = HF'F'_m
\eeq
as claimed.

(ii) Suppose that $\mathrm{Cl}_{\ab(m)}(H)$ is finitely generated. Then $A$ is finite by Theorem \ref{clvfin}(ii).

Conversely, assume that $A$ is finite. By Theorem \ref{clvfin}(iii), it suffices to show that $\mathrm{Cl}_{\ab(m)}(H)$ has finite index in $F$. 
It follows from (\ref{proabm1}) and the homomorphism theorems that
\beq
\label{proabm2}
F/\mathrm{Cl}_{\ab(m)}(H) = F/\pi\inv(\pi(H)) \cong (F/F'F^m)/\pi(H).
\eeq
Since $F/F'F^m$ is finite, we are done.

(iii) It follows immediately from part (ii).

(iv) If $\mathrm{Cl}_{\ab(m)}(H)$ is finitely generated, then $A$ is finite by part (ii). Then 
$\mathrm{Cl}_{\ab(m)}(H)$
has finite index in $F$ and
we can effectively compute the finite group $(F/F'F^m)/\pi(H)$, which is isomorphic to $F/\mathrm{Cl}_{\ab(m)}(H)$ by (\ref{proabm2}). Now we proceed as in the proof of Theorem \ref{proab}(v).

(v) By Theorem \ref{clvfin}(iv).

(vi) Suppose that $H$ is $\mathbf{Ab}(m)$-dense. Let $p$ be any prime dividing $m$. Since $\mathbf{Ab}(p)\subseteq \mathbf{Ab}(m)$, then $H$ is $\mathbf{Ab}(p)$-dense.

Suppose now that $H$ is not $\mathbf{Ab}(m)$-dense.  Then there exists a normal subgroup $N$ of $F$ such that $F/N\in \mathbf{Ab}(m)$ and $HN < F$. Note that $F/N$ is a finite abelian group and $HN$ is contained in a maximal subgroup $M$ of $F$. Since $F/N$ is abelian, we have $F' \leq N \leq M$ and so $M \lhd F$.
Now $M/N$ is a maximal subgroup of $F/N \in \mathbf{Ab}(m)$, hence $[F/N : M/N] = p$ for some prime $p$ dividing $m$. Thus $[F : M] = p$ and so $F/M \cong C_p \in  \mathbf{Ab}(p)$. However $HM=M <  F$ and therefore $H$ is not $ \mathbf{Ab}(p)$-dense.

(vii) By the argument in the proof of Theorem \ref{proab}(vii).

(viii) Since $\cl_{\ab(m)}(H) = \pi\inv(\pi(H))$ by (\ref{proabm1}), it suffices to note that the membership problem is decidable for every finitely generated subgroup of the group
$$F/(F'F^m) \cong \bigoplus_{a\in A} C_m,$$ 
which is straightforward.

(ix) By part (viii) and Proposition \ref{dmp}.
\qed

We can now relate the pro-$\ab$ topology with the pro-$\ab(p^k)$ topologies, where $p$ is a prime and $k$ is a positive integer:

\bp
\label{p:abelem} 
Let $H \leq_{f.g.} F_A$. Then:
\bi
\item[(i)] $H$ is $\mathbf{Ab}$-clopen  if and only if there exist $r\geq 1$,  primes $p_1,\dots,p_r$, positive integers $k_1,\dots,k_r$ and $\ab(p_i^{k_i})$-clopen subgroups $H_i$ of $F_A$ for $1\leq i \leq r$ such that 
$H=\ds\cap_{i=1}^r H_i$;
\item[(ii)] $H$ is  $\mathbf{Ab}$-dense if and only if  $H$ is 
$\ab(p)$-dense  for every prime $p$;
\item[(iii)] $${\rm{Cl}}_{\mathbf{Ab}}(H)=\bigcap_{p\in \mathbb{P}}\bigcap_{k\in \mathbb{N}}{\rm{Cl}}_{\ab(p^{k})}(H)=\bigcap_{p \in \mathbb{P}}\bigcap_{k\in \mathbb{N}}HF_A'F_A^{p^k}.$$
\ei
\ep

\proof
Write $F = F_A$. 

(i) Suppose that there exist $r\geq 1$,  primes $p_1,\dots,p_r$, positive integers $k_1,\dots,k_r$ and $\ab(p_i^{k_i})$-clopen subgroups $H_i$ of $F$ for $1\leq i \leq r$ such that 
$H= \cap_{i=1}^r H_i$. Let $1\leq i \leq r$. Then $F/\textrm{Core}_F(H_i) \in \ab(p_i^{k_i}) \subset \ab$ and so 
$H_i$ is $\mathbf{Ab}$-clopen. Hence $H=\cap_{i=1}^rH_i$ is $\mathbf{Ab}$-clopen.

Suppose now that $H$ is $\mathbf{Ab}$-clopen. 
If $H=F$ the result is trivial: take $r=1$, $H_1=H$ and $p_1$ any prime. Suppose $H\neq F$.  Let $M(H)=F/\textrm{Core}_F(H)$ and  $\varphi: F\rightarrow M(H)$ be the canonical surjective homomorphism. Note that $H=\p^{-1}(\p(H))$.

Also, since $H$ is $\mathbf{Ab}$-clopen, $M(H)$ is a (nontrivial) finite abelian group, hence $F' \unlhd \textrm{Core}_F(H) \leq H$ and so $H\lhd F$ and $\textrm{Core}_F(H)=H$. In particular $\p(H)$ is trivial.

 Since  $M(H)$ is a nontrivial finite abelian group, there exist $r\geq 1$,  primes $p_1,\dots, p_r$ and positive integers $k_1,\dots,k_r$, $a_1,\dots,a_r$ such that $$M(H)=\times_{i=1}^r C_{p_i^{k_i}}^{a_i}.$$ For $1\leq i \leq r$, let $G_i=C_{p_i^{k_i}}^{a_i}$ and $\gamma_i:M(H)\rightarrow G_i$ be the projection of $M(H)$ onto its $i$-th component.
Since $\cap_{i=1}^r\gamma_i^{-1}(1) = 1 = \p(H)$, we have 
$$ H =\p^{-1}(\p(H))=\p^{-1}\left( \ds\bigcap_{i=1}^r\gamma_i^{-1}(1)\right) =\bigcap_{i=1}^r \p^{-1}(\gamma_i^{-1}(1)).$$
Since $\gamma_i^{-1}(1)= \textrm{Ker} \ \gamma_i \unlhd M(H)$, we get $\p^{-1}(\gamma_i^{-1}(1)) \unlhd F$. 
Finally, $F/\p^{-1}(\gamma_i^{-1}(1)) \cong G_i \in \ab(p_i^{k_i})$ and we obtain that $\p^{-1}(\gamma_i^{-1}(1))$ is $\ab(p_i^{k_i})$-clopen. Setting $H_i=\p^{-1}(\gamma_i^{-1}(1))$, we get the result. 

(ii) Let $p$ be a prime. Since  $\ab(p) \subseteq \textbf{Ab}$, if $H$ is $\textbf{Ab}$-dense  then $H$ is $\ab(p)$-dense. 

We now prove the other direction. 
Suppose that $H$ is not $\textbf{Ab}$-dense. Then there is a proper $\textbf{Ab}$-clopen subgroup $K$ of $F$ containing $H$. By part (i) there exist a positive integer $r$, primes $p_1,\dots,p_r$, positive integers $k_1,\dots,k_r$ and subgroups $K_1,\dots,K_r$ of $F$ such that, for $1\leq i \leq r$, $K_i$ is $\ab(p_i^{k_i})$-clopen and $K=\cap_{i=1}^r K_i$. As $K$ is proper, there exists $1\leq i\leq r$ such that $K_i$ is proper.  Since $H$ is contained in $K_i$, $H$ is not $\ab(p_i^{k_i})$-dense. Therefore $H$ is not $\ab(p_i)$-dense by Theorem \ref{proabm}(vi).
 
(iii) First note that as for every prime $p$ and every positive integer $k$, we have $\la C_{p^{k}} \ra \subseteq \mathbf{Ab}$, we obtain 
 $${\rm{Cl}}_{\mathbf{Ab}}(H)\leq \bigcap_{p\in \mathbb{P}}  \bigcap_{k\in \mathbb{N}}{\rm{Cl}}_{\ab(p^{k})}(H).$$ We consider the reverse inclusion. 
By Theorem \ref{sclosed}, ${\rm{Cl}}_{\mathbf{Ab}}(H)$ is the intersection of all the $\ab$-(cl)open subgroups of $F$ containing $H$. Let $K$ be an $\textbf{Ab}$-clopen subgroup of $F$ containing $H$. By part (i) there exist $r\geq 1$, primes $p_1,\dots,p_r$, positive integers $k_1,\dots,k_r$ and $\ab(p_i^{k_i})$-clopen subgroups $K_i$ of $F$ for $1\leq i \leq r$ such that 
$K=\cap_{i=1}^r K_i$.  In particular, for each $1\leq i \leq r$,  ${\rm{Cl}}_{\ab(p_i^{k_i})}(H)\leq K_i$. Hence $$\bigcap_{p\in \mathbb{P}}\bigcap_{k\in \mathbb{N}}{\rm{Cl}}_{\ab(p^{k})}(H)\leq K$$
for every $\textbf{Ab}$-clopen subgroup $K$ of $F$ containing $H$, yielding
$$\bigcap_{p\in \mathbb{P}}\bigcap_{k\in \mathbb{N}}{\rm{Cl}}_{\ab(p^{k})}(H)
\leq{\rm{Cl}}_{\mathbf{Ab}}(H).$$
The second equality follows from Theorem \ref{proabm}(i). 
\qed

Let ${\bf N}$ denote the pseudovariety of all finite nilpotent groups. For every prime $p$, let ${\bf G}_p$ denote the pseudovariety of all finite $p$-groups.
Using results from \cite{MSW}, we obtain:

\bc
\label{nab}
Let $H \leq_{f.g.} F_A$. Then the following conditions are equivalent:
\bi
\item[(i)] $H$ is {\bf N}-dense;
\item[(ii)] $H$ is ${\bf G}_p$-dense for every prime $p$;
\item[(iii)] $H$ is {\bf Ab}$(p)$-dense for every prime $p$;
\item[(iv)] $H$ is {\bf Ab}-dense.
\ei
\ec

\proof
Suppose that $A$ is infinite. By the proof of Theorem \ref{proab}(vii) (which adapts directly to prove Theorem \ref{proabm}(vii)), $H$ is not $\ab(p)$-dense for any prime $p$. Since $\ab(p)$ is contained in {\bf N}, ${\bf G}_p$ and {\bf Ab}, it follows that $H$ is neither {\bf N}-dense nor ${\bf G}_p$-dense nor {\bf Ab}-dense.

Thus we may assume that $A$ is finite. Now
(i) $\iff$ (ii) follows from \cite[Corollary 4.2(1)]{MSW},
(ii) $\iff$ (iii) follows from \cite[Corollary 3.5]{MSW} and
(iii) $\iff$ (iv) follows from Proposition \ref{p:abelem}(ii).
\qed

\section*{Acknowledgements}

The first author acknowledges support from the Centre of Mathematics of the University of Porto, which is financed by national funds through the Funda\c c\~ao para a Ci\^encia e a Tecnologia, I.P., under the project with references UIDB/00144/2020 and  UIDP/00144/2020.

The second author acknowledges support from the Centre of Mathematics of the University of Porto, which is financed by national funds through the Funda\c c\~ao para a Ci\^encia e a Tecnologia, I.P., under the project with reference UIDB/00144/2020. 

The third author was supported by the Engineering and Physical Sciences Research Council, grant number 
EP/T017619/1.

\vspace{1cm}

{\sc Claude Marion, Centro de
Matem\'{a}tica, Faculdade de Ci\^{e}ncias, Universidade do
Porto, R. Campo Alegre 687, 4169-007 Porto, Portugal}

{\em E-mail address}: claude.marion@fc.up.pt

\bigskip

{\sc Pedro V. Silva, Centro de
Matem\'{a}tica, Faculdade de Ci\^{e}ncias, Universidade do
Porto, R. Campo Alegre 687, 4169-007 Porto, Portugal}

{\em E-mail address}: pvsilva@fc.up.pt

\bigskip

{\sc Gareth Tracey, Mathematics Institute, University of Warwick, Coventry CV4 7AL, U.K.}

{\em E-mail address}: Gareth.Tracey@warwick.ac.uk

\end{document}